\documentclass[10pt,preprint]{amsart}
\usepackage{amssymb,amsfonts,amsthm,amsmath,amscd,relsize}
\usepackage{hyperref}
\usepackage{enumerate}
\usepackage{comment}

 \usepackage{epsfig}
 \usepackage{graphicx}
 \usepackage{lipsum}
\usepackage{xcolor}
\usepackage{float}
\usepackage{soul}
\usepackage{adjustbox}

\usepackage{enumitem}

\usepackage{pgf,tikz,pgfplots}
\pgfplotsset{compat=1.13}
\usepackage{mathrsfs}
\usetikzlibrary{arrows}
\usetikzlibrary{calc}

\definecolor{qqqqff}{rgb}{0,0,1}
\definecolor{qqwuqq}{rgb}{0,0.39215686274509803,0}
\definecolor{uuuuuu}{rgb}{0.266,0.266,0.266}

\def\R{\mathbb R}

\def\A{\textbf{\textit{A}}}
\def\B{\textbf{\textit{B}}}
\def\C{\textbf{\textit{C}}}
\def\D{\textbf{\textit{D}}}
\def\T{\textbf{\textit{T}}}

\theoremstyle{plain}
\newtheorem*{theorem}{Theorem}
\newtheorem*{corollary}{Corollary}

\theoremstyle{definition}

\theoremstyle{remark}

\makeatletter
\@addtoreset{footnote}{page}
\makeatother

\definecolor{ttzzqq}{rgb}{0.2,0.6,0}

\begin{document}


\title{A simple sum for simplices}

\author{Christian Aebi and Grant Cairns}

\address{Coll\`ege Calvin, Geneva, Switzerland 1211}
\email{christian.aebi@edu.ge.ch}
\address{Department of Mathematics, La Trobe University, Melbourne, Australia 3086}
\email{G.Cairns@latrobe.edu.au}

\begin{abstract} We give a  vector identity for $n+2$ points in $\R^n$. It follows as a corollary that 
when $n$ is odd the sum of the signed volumes of the $n$-simplices is zero, and when $n$ is even, the alternating sum of the signed volumes is zero.
\end{abstract}

\maketitle

\section{Introduction}

Consider a planar quadrilateral with vertices having position vectors $\A,\B,\C,\D$, in cyclic order, as in the diagram. 

\begin{figure}[htp]
\begin{center}
\includegraphics[scale=.75]{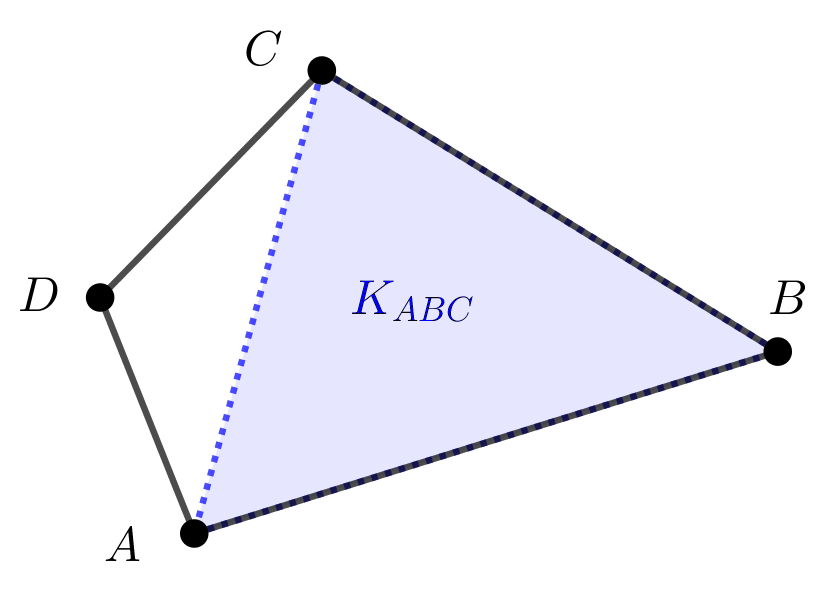}
\caption{Illustration of $K_{ABC}$ of the $ABCD$ quadrilateral}\label{F:quad}
\end{center}
\end{figure}

Let $K_{ABC}$ denote the area of the triangle $ABC$, and use similar notation for the other triangles. The identity 
\begin{equation}\label{E:college}
K_{BCD}\, \A -K_{ACD}\, \B+K_{ABD}\, \C  -K_{ABC}\, \D=0
\end{equation}
 was given in \cite{CA} where it was proved as a consequence of the Jacobi vector triple product identity in $\R^3$.

One has $2K_{ABC}=\det[\B-\A,\C-\A]$, and similarly for the other triangles, so \eqref{E:college} gives
\begin{equation}\label{E:collegedet}
\begin{gathered}
\det[\C-\B,\D-\B]\, \A -\det[\D-\C,\A-\C]\, \B\\\quad\,+\det[\A-\D,\B-\D]\, \C  -\det[\B-\A,\C-\A]\, \D=0.
\end{gathered}
\end{equation}
The object of this paper is to generalize this fact to arbitrary dimension $n$. Consider points $\A_0,\dots,\A_{n+1}$ in $\R^n$ and think of them as column vectors. For $i=0,\dots,n+1$, consider the 
$n\times n$ matrix 
\[
M_i=\big[\A_{i+2}-\A_{i+1}\, \big\vert\,  \A_{i+3}-\A_{i+1}\, \big\vert\, \dots \,\big\vert\, \A_{i+n+1}-\A_{i+1}\big],
\]
where the indices are computed modulo $n+2$, and let $\Delta_i=\det M_i$. Here below is the main result of this paper.

\begin{theorem} In the above notation, one has
\begin{equation}\label{Id1}
\sum_{i=0}^{n+1}\,  (-1)^{i(n+1)} \, \Delta_i  \, \A_i=0.
\end{equation}
\end{theorem}

\begin{figure}[htp]
\begin{center}
\includegraphics[scale=.75]{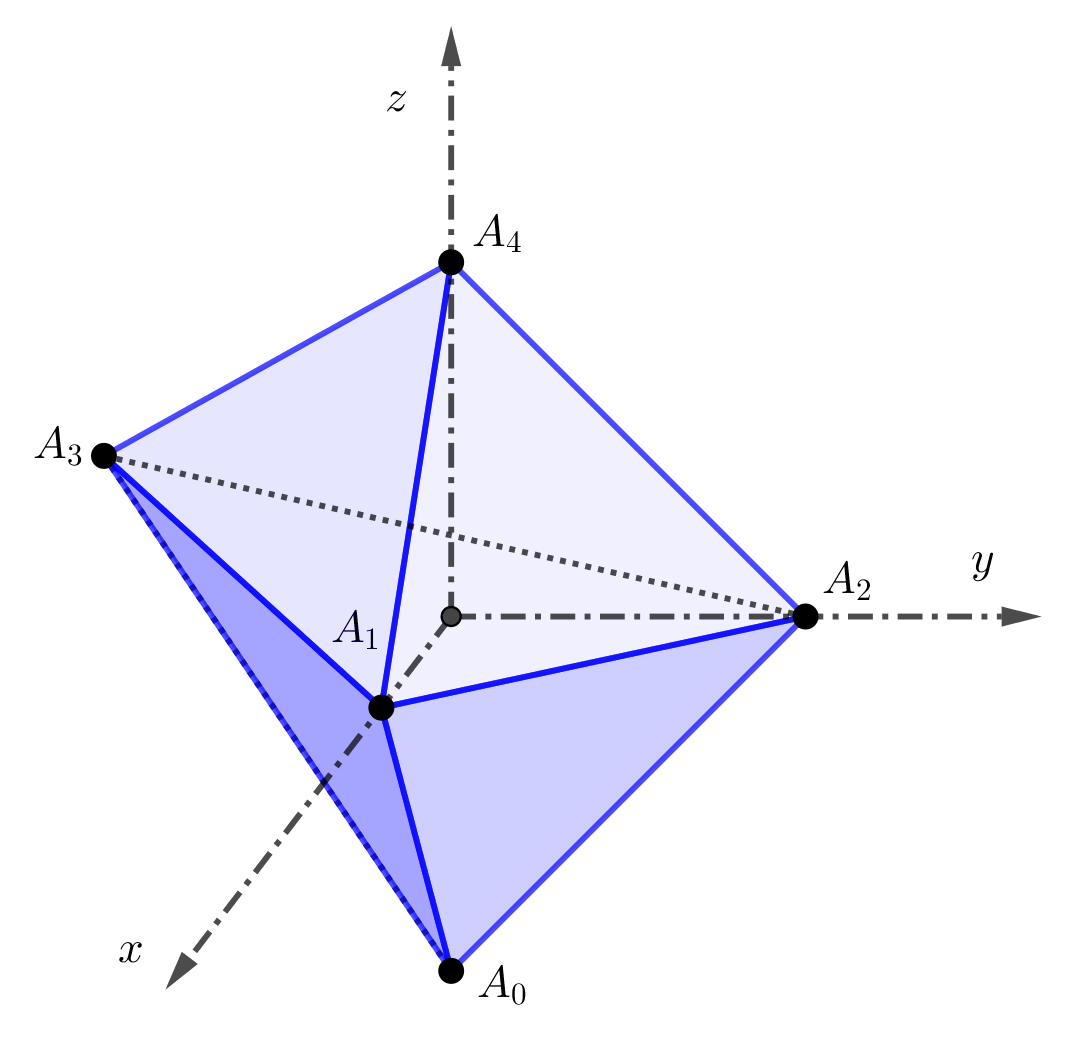}
\caption{Applying (\ref{Id1}) on a triangular bipyramid }\label{F:bipyr}
\end{center}
\end{figure}

Before proving this result, let us consider some special cases. For $n=1$, the above theorem gives
\[
(\A_2-\A_1)\, \A_0+(\A_0-\A_2)\, \A_1+(\A_1-\A_0)\, \A_2=0,
\]
which is obvious. 
For $n=2$, the theorem gives us \eqref{E:collegedet}.
For $n=3$, we have  5 points in $\R^3$, which we may regard as the vertices of a (possibly degenerate) polyhedron, and the theorem gives
 \begin{align*}
\det&[\A_2-\A_1,\A_3-\A_1,\A_4-\A_1]\, \A_0+\det[\A_3-\A_2,\A_4-\A_2,\A_0-\A_2]\, \A_1\\
&+\det[\A_4-\A_3,\A_0-\A_3,\A_1-\A_3]\, \A_2+\det[\A_0-\A_4,\A_1-\A_4,\A_2-\A_4]\, \A_3\\
&+\det[\A_1-\A_0,\A_2-\A_0,\A_3-\A_0]\, \A_4=0.
  \end{align*}
Here for each $i$, the coefficient of $ \A_i$ is $6$ times the signed volume of the tetrahedron defined by the other four vertices.
For example, in Figure~\ref{F:bipyr}, for $\A_0=(0, 0, -1), \A_1=(1, 0, 0),  \A_2=(0, 1,0), \A_3=(-1, -1, 0), \A_4=(0, 0,1)$, we obtain
\[
3\A_0-2 \A_1-2\A_2-2 \A_3+3\A_4=0.
 \]

In dimension $n$, the convex hull of $n+1$ points is a (possibly degenerate) \emph{$n$-simplex}. The coefficient $\Delta_i$  in \eqref{Id1} is the signed volume of the $n$-simplex defined by the points other than $\A_i$. In particular, $\Delta_i$ is unchanged by translation. So  translating by a nonzero vector $\T$,  \eqref{Id1} gives 
$\sum_{i=0}^{n+1}\,  (-1)^{i(n+1)} \, \Delta_i  (A_i+\T)=0$. Then subtracting \eqref{Id1}  and taking the coefficient of $\T$ gives the following scalar identity.

\begin{corollary} $\displaystyle
\sum_{i=0}^{n+1}\,  (-1)^{i(n+1)} \, \Delta_i  =0$.
\end{corollary}

In other words, given $n+2$ points in $\R^n$, when $n$ is odd the sum of the signed volumes of the $n$-simplices is zero, and when $n$ is even, the alternating sum of the signed volumes is zero.

\section{Multilinear algebra}

Our  proof of the theorem is a simple argument using multilinear algebra. Let us summarize the well known  basic ideas  we require.
Consider real vector spaces $V$ and $W$. Suppose that $k$ is a positive integer. Recall that a function $f:V^{k}\to W$ is \emph{multilinear} if it is linear in each variable with the other variables held constant; for a gentle introduction, see \cite[Chap.~3]{Tu}. 
A multilinear function $f:V^{k}\to W$ is \emph{alternating} if for all elements $\A_0,\A_1,\dots,\A_{k-1}\in V$ and all permutations $\sigma$ of $\{0,1,\dots,k-1\}$, one has
\[
f(\A_{\sigma(0)},\A_{\sigma(1)},\dots,\A_{\sigma(k-1)})=\text{Sgn}(\sigma) f(\A_0,\A_1,\dots,\A_{k-1}),
\]
where $\text{Sgn}(\sigma)$ denotes the sign of $\sigma$.
For example, the determinant function $\det: V^{n}\to \R, (\A_0,\dots,\A_{n-1})\mapsto\det[\A_0\, \big\vert\, \dots\, \big\vert\, \A_{n-1}]$ is an alternating multilinear function of the column vectors; see \cite[Chap.~XIII]{Lang}.

There are several known sets of generators for the symmetric group $S_{k}$ of permutations of $\{0,1,\dots,k-1\}$; a nice survey is given in \cite{Conrad}. In particular, $S_{k}$ is generated by the cycle $(0,1,\dots,k-1)$ and the transposition $(0,1)$; see \cite[Theorem~2.5]{Conrad}. Thus, in order to show that a multilinear function $f$ is alternating it suffices to show that for all  $\A_0,\A_1,\dots,\A_{k-1}\in V$,
\begin{enumerate}
\item $f(\A_1,\A_2,\dots,\A_{k-1},\A_{0})=(-1)^{k-1}f(\A_0,\A_1,\dots,\A_{k-1})$,
\item $f(\A_1,\A_0,\A_2,\dots,\A_{k-1})=-f(\A_0,\A_1,\A_2,\dots,\A_{k-1})$.
\end{enumerate}
Note that from (b) we have
\begin{enumerate}
\item[(b$'$)] $f(\A_0,\A_0,\A_2,\A_3,\dots,\A_{k-1})=0$.
\end{enumerate}
Conversely,  it is easy to see that if (b$'$) holds for all  $\A_0,\A_2,\dots,\A_{k-1}\in V$, then (b) follows.
So in order to show a multilinear function $f$ is alternating it suffices to verify conditions (a) and (b$'$). Note that it follows that if $f$ is alternating, and if $\A_i=\A_j$ for some $i\not=j$, then $f(\A_0,\A_1,\dots,\A_{k-1})=0$; indeed, one can just permute $\A_i,\A_j$ to the extreme left, and employ (b$'$).

Finally, a key fact about alternating multilinear functions that we will use below is that if $k>\dim V$, then $f$ is identically zero. In the literature, this fact can be quickly deduced once one has constructed the exterior algebra on $V$, but we don't do that here as we won't require the exterior product. Instead, one can use the following straight forward proof. First choose a basis for $V$. Using multilinearity, the image of alternating  function $f$ is determined by its values on the basis elements. But if $\A_0,\A_1,\dots,\A_{k-1}$ are basis elements and $k>\dim V$, then by the pigeonhole principle, there must be a repetition of one of the basis elements. It follows that as $f$ is alternating, $f(\A_0,\A_1,\dots,\A_{k-1})=0$.

\section{Proof of the Theorem}
Let $n$ be an arbitrary positive integer, let $V=\R^n$ and consider the function $f:V^{n+2}\to V$ defined by
\[
f(\A_0,\A_1,\dots,\A_{n+1})=\sum_{i=0}^{n+1} (-1)^{i(n+1)} \Delta_i  \, \A_i.
\]
Because $\det$ is multilinear, and because for each $i$, the variable $\A_i$ does not occur in $\Delta_i$, it follows that $f$ is multilinear. We will show that $f$ is alternating,  and hence identically zero. Condition (a) is immediate from the definition of $f$. So it remains to prove (b$'$).
 Using the fact that the determinant  is an alternating multilinear function of the column vectors and computing the indices modulo $n+2$, we have
\begin{align*}
 \Delta_i=& \det\big[ \A_{i+2}-\A_{i+1}\, \big\vert\, \A_{i+3}-\A_{i+1}\, \big\vert\, \dots \,\big\vert\, \A_{i+n+1}-\A_{i+1}\big]
\\
=&\det\big[\A_{i+2}\, \big\vert\, \A_{i+3}\, \big\vert\, \dots \,\big\vert\,\A_{i+n+1}\big]\\
 &-\sum_{j=2}^{n+1}\det\big[\A_{i+2}\, \big\vert\, \dots\,\big\vert\, \A_{i+j-1}\, \big\vert\, \A_{i+1} \,\big\vert\, \A_{i+j+1}\, \big\vert\, \dots \,\big\vert\,\A_{i+n+1}\big].
  \end{align*}
Moving the column $\A_{i+1}$ in the summation $j-2$ positions to the far left, we have
 \begin{align*}
 \Delta_i=&\det\big[\A_{i+2}\, \big\vert\, \A_{i+3}\, \big\vert\, \dots \,\big\vert\,\A_{i+n+1}\big]\\
 &-\sum_{j=2}^{n+1}(-1)^{j-2}\det\big[\A_{i+1}\, \big\vert\, \dots\,\big\vert\, \A_{i+j-1}\, \big\vert\, \widehat{A_{i+j}} \,\big\vert\, \A_{i+j+1}\, \big\vert\, \dots \,\big\vert\,\A_{i+n+1}\big]\\
 =&\sum_{j=1}^{n+1}(-1)^{j-1}\det\big[\A_{i+1}\, \big\vert\, \dots\,\big\vert\, \A_{i+j-1}\, \big\vert\, \widehat{A_{i+j}} \,\big\vert\, \A_{i+j+1}\, \big\vert\, \dots \,\big\vert\,\A_{i+n+1}\big],
 \end{align*}
 where in the above, the hat symbol indicates that the term has been omitted. In particular,
 
 \begin{align*}
 \Delta_0=&\sum_{j=1}^{n+1}(-1)^{j-1}\det\big[  \A_{1}\, \big\vert\, \dots\,\big\vert\, \widehat{\A_{j}} \, \big\vert\, \dots \,\big\vert\,\A_{n+1}\big].
 \end{align*}
and
\begin{align*}
 \Delta_1=&\sum_{j=1}^{n+1}(-1)^{j-1}\det\big[  \A_{2}\, \big\vert\, \dots\,\big\vert\, \widehat{\A_{j+1}} \, \big\vert\, \dots \,\big\vert\,\A_{n+1}\,\big\vert\,\A_{0}\big].
 \end{align*}
 Now suppose $\A_1=\A_0$. Moving $\A_0$ to the far left and replacing it by $\A_1$, and then adjusting $j$, we have
 \begin{align*}
 \Delta_1=&\sum_{j=0}^{n}(-1)^{n+j}\det\big[  \A_{1}\, \big\vert\, \dots\,\big\vert\, \widehat{\A_{j+1}} \, \big\vert\, \dots \,\big\vert\,\A_{n+1}\big]\\
 =&\sum_{j=1}^{n+1}(-1)^{n+j-1}\det\big[  \A_{1}\, \big\vert\, \dots\,\big\vert\, \widehat{\A_{j}} \, \big\vert\, \dots \,\big\vert\,\A_{n+1}\big].
 \end{align*}
So $\Delta_0+ (-1)^{n+1}\Delta_1 =0$. Since $\A_1=\A_0$, we have $\Delta_i=0$ for all $i\ge 2$. Hence, from the definition of $f$,
\[
f(\A_0,\A_0,\A_2,\dots,\A_{n+1})=(\Delta_0+ (-1)^{n+1}\Delta_1)\A_{0}=0,
\]
as required.

\bibliographystyle{amsplain}
{}


\end{document}